\documentclass[10pt]{article}
\usepackage{amssymb}
\usepackage{amsmath}
\usepackage{amsthm}
\usepackage{mathabx} 
\usepackage[backref=page,colorlinks=true]{hyperref} 
\usepackage[capitalise]{cleveref} 

\makeatletter

\def\myMRbibitem{\@ifnextchar[\my@lbibitem\my@bibitem}

\def\mybiblabel#1#2{\@biblabel{{\hyperref{http://www.ams.org/mathscinet-getitem?mr=#1}{}{}{#2}}}}

\def\myhyperanchor#1{\Hy@raisedlink{\hyper@anchorstart{cite.#1}\hyper@anchorend}}

\def\my@lbibitem[#1]#2#3#4\par{%
    \item[\mybiblabel{#2}{#1}\myhyperanchor{#3}\hfill]#4%
    \@ifundefined{ifbackrefparscan}{}{\BR@backref{#3}}%
    \if@filesw{\let\protect\noexpand\immediate
       \write\@auxout{\string\bibcite{#3}{#1}}}\fi\ignorespaces%
}

\def\my@bibitem#1#2#3\par{%
    \refstepcounter\@listctr
    \item[\mybiblabel{#1}{\the\value\@listctr}\myhyperanchor{#2}\hfill]#3%
    \@ifundefined{ifbackrefparscan}{}{\BR@backref{#2}}%
    \if@filesw\immediate\write\@auxout
        {\string\bibcite{#2}{\the\value\@listctr}}\fi\ignorespaces%
}

\makeatother

\theoremstyle{plain}
	
	\crefname{mainthm}{Theorem}{Theorems} 
	\newtheorem{theorem}{Theorem}[section]
	
	\newtheorem{corollary}[theorem]{Corollary}
	\newtheorem{proposition}[theorem]{Proposition}

\theoremstyle{remark}
	\newtheorem{remark}[theorem]{Remark}
	
	\newtheorem*{ack}{Acknowledgement}

\newcommand{\closeremark}{\hfill$\triangleleft$}  

\crefname{subsection}{\S}{\S\S} 
\Crefname{subsection}{\S}{\S\S} 

\numberwithin{equation}{section}

\renewcommand{\epsilon}{\varepsilon}
\renewcommand{\phi}{\varphi}

\newcommand{\cA}{\mathcal{A}}\newcommand{\cB}{\mathcal{B}}\newcommand{\cC}{\mathcal{C}}

\newcommand{\cI}{\mathcal{I}}

\newcommand{\cP}{\mathcal{P}}\newcommand{\cQ}{\mathcal{Q}}\newcommand{\cR}{\mathcal{R}}
\newcommand{\cU}{\mathcal{U}}

\newcommand{\R}{\mathbb R}

\newcommand{\Z}{\mathbb Z}

\newcommand{\SL}{\mathrm{SL}}

\newcommand{\m}{\mathfrak{m}}    
\newcommand{\id}{\mathrm{id}}  

\DeclareMathOperator\Aut{Aut}

\newcommand{\bysame}{\makebox[3em]{\hrulefill}\thinspace. }

\newcommand{\tribar}[1]{\mathopen{| {\kern -1.5pt} | {\kern -1.5pt} |} {#1} \mathclose{| {\kern -1.5pt} | {\kern -1.5pt} |}}
\newcommand{\bigtribar}[1]{\mathopen{\big|{\kern -1.5pt}\big|{\kern -1.5pt}\big|}{#1}\mathclose{\big|{\kern -1.5pt}\big|{\kern -1.5pt}\big|}}
\newcommand{\Bigtribar}[1]{\mathopen{\Big|{\kern -1.5pt}\Big|{\kern -1.5pt}\Big|}{#1}\mathclose{\Big|{\kern -1.5pt}\Big|{\kern -1.5pt}\Big|}}
\newcommand{\biggtribar}[1]{\mathopen{\bigg|{\kern -1.5pt}\bigg|{\kern -1.5pt}\bigg|}{#1}\mathclose{\bigg|{\kern -1.5pt}\bigg|{\kern -1.5pt}\bigg|}}
\newcommand{\biangle}[1]{\mathopen{\langle {\kern -1.5pt} \langle}{#1} \mathclose{\rangle {\kern -1.5pt} \rangle}}

\newcommand{\arxiv}[1]{Preprint \href{http://arxiv.org/abs/#1}{arXiv:{#1}}}

\renewcommand*{\backref}[1]{}
\renewcommand*{\backrefalt}[4]{\quad \tiny 
    \ifcase #1 (Not cited.)%
    \or        (Cited on page~#2.)%
    \else      (Cited on pages~#2.)%
    \fi}

\begin{document}

\author{Jairo Bochi\footnote{Partially supported by CNPq (Brazil) and FAPERJ (Brazil).} \quad \& \quad Andr\'es Navas\footnote{Partially supported by the Fondecyt Project 1120131 and the  ``Center of Dynamical Systems and 
Related Topics" (DySyRF; ACT-Project 1103, Conicyt).}}

\title{Almost Reduction and Perturbation of Matrix Cocycles}

\date{August 14, 2013}

\maketitle
\thispagestyle{empty} 

\begin{abstract}
In this note, we show that if all Lyapunov exponents 
of a matrix cocycle vanish, then it can be perturbed to become 
cohomologous to a cocycle taking values in the orthogonal group.
This extends a result of Avila, Bochi and Damanik to general 
base dynamics and arbitrary dimension.
We actually prove a fibered version of this result, and 
apply it to study the existence of dominated splittings into 
conformal subbundles for general matrix cocycles.
\end{abstract}

\section{From zero Lyapunov exponents to rotation cocycles} 

\subsection{Basic definitions}

Let $F \colon \Omega \to \Omega$ be a homeomorphism of a compact metric space $\Omega$.
Let $V$ be a finite-dimensional real vector bundle over $\Omega$, whose fiber over $\omega$ 
is denoted by $V_\omega$. 
Let $\mathcal{A}$ be a vector-bundle automorphism that fibers over $F$; this means 
that the restriction of $\cA$ to each fiber $V_\omega$ 
is a linear automorphism $A(\omega)$ onto $V_{F\omega}$.
In the case of trivial vector bundles, $\cA$ is usually called 
a \emph{linear cocycle}.

As a convention, automorphisms of $V$ will be denoted by calligraphic letters,
and the restrictions to the fibers will be denoted by the corresponding roman letters.
Analogously, for any integer $n$,  the restriction of the power $\cA^n$ to the fiber $V_\omega$
is denoted by $A^n(\omega)$; 
thus $A^n(\omega) = A(F^{n-1}\omega) \circ \cdots \circ A(\omega)$ for $n >0$.

A \emph{Riemannian metric} on $V$ is a continuous choice of inner product 
$\langle \mathord{\cdot}, \mathord{\cdot} \rangle_\omega$ 
on each fiber $V_\omega$.
It induces a \emph{Riemannian norm}
$\|v\|_\omega = \sqrt{\langle v, v \rangle_\omega}${\,}.
Given a linear map $L \colon V_{\omega} \to V_{\omega'}$,
its \emph{norm} $\|L\|$ and its \emph{mininorm} $\m(L)$
are defined respectively as the supremum and the infimum of $\| Lv\|_{\omega'}$
over all unit vectors $v \in V_\omega$.

Let $\Aut(V,F)$ denote the space of all automorphisms of $V$ that fiber over $F$,
endowed with the topology induced by the distance 
$d(\cA,\cB) = \sup_\omega \|A(\omega) - B(\omega)\|$,
for some choice of a Riemannian norm on $V$.

\subsection{Uniform subexponential growth and its consequences}

Define
$$
\lambda^+(\cA) = \lim_{n \to +\infty} \frac{1}{n}  \sup_{\omega \in \Omega} \log \|A^n (\omega) \| 
\quad \text{and} \quad
\lambda^-(\cA) = \lim_{n \to +\infty} \frac{1}{n}  \inf_{\omega \in \Omega} \log \m(A^n (\omega)) \, .
$$
which exist by subadditivity and supraadditivity, respectively.

If $\mu$ is ergodic probability measure for $F \colon\Omega\to \Omega$, 
then there are constants $\lambda^+(\cA,\mu)$, $\lambda^-(\cA,\mu)$,
called the \emph{top} and \emph{bottom Lyapunov exponents},
such that, for $\mu$-almost every $\omega \in \Omega$,
$$
\frac{1}{n} \log \| A^n(\omega) \| \to \lambda^+(\cA,\mu) 
\quad \text{and} \quad
\frac{1}{n} \log \m(A^n(\omega))   \to \lambda^-(\cA,\mu) 
\quad \text{as $n\to+\infty$.}
$$
Moreover, the following ``variational principle'' 
holds\footnote{This follows from \cite[Thrm.~1]{Schreiber} or \cite[Thrm.~1.7]{SS}. 
Although these references assume $\Omega$ to be compact metrizable, the proofs also 
work for compact Hausdorff $\Omega$. (See also the proof of Proposition~1 in 
\cite{AB}.) A particular case was considered in \cite{Furman}.}:
\begin{equation}\label{e.variational}
\lambda^+(\cA) = \sup_\mu \lambda^+(\cA,\mu) \quad \text{and} \quad
\lambda^-(\cA) = \inf_\mu \lambda^-(\cA,\mu) \, .
\end{equation}
where $\mu$ runs over all invariant ergodic probabilities for $F$.

\medskip

Let us say that the automorphism $\cA$ has \emph{uniform subexponential growth} if
$\lambda^+(\cA) = \lambda^-(\cA) = 0$.
By (\ref{e.variational}), this is equivalent to the vanishing of all Lyapunov exponents with respect 
to all ergodic probability measures. 

\medskip

Our first result is:

\begin{theorem}\label{t.zero}
Assume that $\cA \in \Aut(V,F)$ has uniform subexponential growth.
Then:
\begin{enumerate}
\item \label{i.z1}
For any $\epsilon>0$,
there exists a Riemannian norm $\tribar{\mathord{\cdot}}$ on $V$ such that
\begin{equation}\label{e.almost_inv_metric} 
e^{-\epsilon} \tribar{v}_\omega < \tribar{A(\omega) v }_{F \omega} < e^\epsilon \tribar{v}_\omega \, , \quad
\text{for all $\omega \in \Omega$, $v\in V_\omega$.}
\end{equation}

\item \label{i.z2}
There exists an arbitrarily small perturbation of $\cA$ that preserves some Riemannian norm on $V$.
\end{enumerate}

\end{theorem}

As we will see, part~(\ref{i.z1}) follows from
a standard construction in Pesin theory, 
and part~(\ref{i.z2}) follows form part~(\ref{i.z1}).
However, the latter implication is \emph{not} straightforward,
because if $\epsilon$ is small then the Riemannian norm constructed in part~(\ref{i.z1})
may be very distorted with respect to a fixed reference Riemannian norm on $V$.

\medskip

For a reformulation of the theorem 
in terms of conjugacy to isometric automorphisms,
see \cref{ss.conjugacy}. 

\medskip

Despite making stringent assumptions about the automorphism $\cA$, \cref{t.zero}
can be used to obtain very strong properties for a dense subset $D$ of $\Aut(V,F)$,
under the assumption that $F$ is uniquely ergodic (or, in some cases, minimal).
More precisely, we show that for every automorphism $\cA$ in the subset $D$ 
there exists a Riemannian metric norm $\tribar{\mathord{\cdot}}$ on $V$ 
and a splitting of $V$ as a Whitney sum of 
$\cA$-invariant subbundles where $\cA$ acts conformally with respect to the norm $\tribar{\mathord{\cdot}}$.  
Moreover, this splitting is either trivial or dominated.
See \cref{ss.subbundles} for details.

\medskip

In the paper \cite{BN}, we prove results 
about cocycles of isometries of spaces of nonpositive curvature
that generalize \cref{t.zero}. Actually, we first obtained \cref{t.zero}
as a corollary of the geometrical results of \cite{BN}. Later, we realized that the constructions
could be modified or adapted to produce an elementary proof of \cref{t.zero}, which 
we present, together with its applications, in this note.

\subsection{Proof of \cref{t.zero}} 

We need a few preliminaries.

Recall that $V$ is a finite-dimensional vector bundle over the compact space $\Omega$.
We choose and fix a Riemannian metric 
$\langle \mathord{\cdot},\mathord{\cdot} \rangle$ on $V$.
Let $\cB$ be an automorphism of $V$ over a homeomorphism $G \colon \Omega \to \Omega$.
The \emph{transpose} of $\cB$ is the automorphism $\cB^*$ over $G^{-1}$ defined by
$$
\langle B(\omega) u, v \rangle_{G\omega} = \langle  u, B^*(G \omega) v \rangle_{\omega}, \quad
\text{for all } u \in V_\omega, \ v \in V_{G\omega}  \, .
$$
If $\cB^* = \cB$ 
(and thus $G$ is the identity), then $\cB$ is called \emph{symmetric}.
An automorphism $\cP$ is called \emph{positive} if it is symmetric and 
$\langle P(\omega) v, v \rangle_\omega > 0$ for all nonzero $v \in V_\omega$.
We write $\cB < \cC$ if $\cB$ and $\cC$ are symmetric and $\cC-\cB$ is positive.

The following \lcnamecref{p.properties} collects some useful properties:

\begin{proposition}\label{p.properties}
\begin{enumerate}
\item \label{i.prop1}
If $\cA$ is any automorphism and $\cB$ is symmetric, then $\cA^* \cB \cA$ is symmetric;
moreover, if $\cB < \cC$, then $\cA^* \cB \cA < \cA^* \cC \cA$. 

\item \label{i.prop2}
Each positive automorphism $\cP$ has a unique positive square root $\cP^{1/2}$;
moreover, $\cP^{1/2}$ commutes with $\cP$, and the map $\cP \mapsto \cP^{1/2}$ is continuous. 
\item \label{i.prop3}
The square root map is monotonic: if $\cP$, $\cQ$ are positive and 
$\cP < \cQ$, then $\cP^{1/2} < \cQ^{1/2}$.
\end{enumerate}
\end{proposition}

Properties (\ref{i.prop1}) and (\ref{i.prop2}) above are easy exercises.
For a proof of property (\ref{i.prop3}), see \cite[p.~9]{Bhatia}.

\begin{proof}[Proof of \cref{t.zero}]
Let $\cA$ be an automorphism of $V$ over the homeomorphism $F$  
having uniform subexponential growth. Fix a small $\epsilon>0$.

To prove part~(\ref{i.z1}),
we will use an standard construction in Pesin theory
called \emph{Lyapunov norms} (see e.g.~\cite[p.~667]{KatokH}). 
Define
\begin{equation}\label{e.Lyapunov_norm}
\tribar{v}^2_\omega := \sum_{n\in \Z} e^{-2\epsilon |n|} \, \| A^n(\omega) v \|^2_{F^n \omega} \, .
\end{equation}
Since the cocycle has uniform subexponential growth,
the series converges uniformly on compact subsets of $V$, 
and hence defines a (continuous) Riemannian norm.
Property~\eqref{e.almost_inv_metric} is straightforward to check.
This proves part~(\ref{i.z1}),

\medskip

To prove part~(\ref{i.z2}), let $\biangle{\mathord{\cdot},\mathord{\cdot}}$ be the inner product
that induces the norm \eqref{e.Lyapunov_norm}.
Then there are positive automorphisms $\cR$, $\cQ$
such that for all $u$, $v\in V_\omega$,
\begin{align}
\biangle{u,v}_\omega                       &= \langle R(\omega)u, v \rangle_\omega \, , \label{e.def_R} \\
\biangle{A(\omega)u,A(\omega)v}_{F \omega} &= \langle Q(\omega)u, v \rangle_\omega \, . \label{e.def_Q}
\end{align}
The almost-invariance property \eqref{e.almost_inv_metric} can now be expressed as:
\begin{equation}\label{e.almost_inv_smart}
e^{-2\epsilon} \cR < \cQ < e^{2\epsilon} \cR \, .
\end{equation}

We want to find an automorphism $\tilde \cA$ over $F$ that is close to $\cA$
and leaves the inner product $\biangle{\mathord{\cdot},\mathord{\cdot}}$ invariant.
As it is straightforward to check, 
invariance means that the automorphism
$\cP = \cA^{-1} \tilde{\cA}$ (over the identity) satisfies:
\begin{equation}\label{e.invariance}
\cP^* \cQ \cP = \cR \, .
\end{equation}
Equivalently,
$$
(\cQ^{1/2} \cP \cQ^{1/2})^*  (\cQ^{1/2} \cP \cQ^{1/2})  = \cQ^{1/2} \cR \cQ^{1/2} \, .
$$
Let us try to find a \emph{positive} solution $\cP$.
Then the relation above becomes
$(\cQ^{1/2} \cP \cQ^{1/2})^2 = \cQ^{1/2} \cR \cQ^{1/2}$,
and using the uniqueness of positive square roots
(property~(\ref{i.prop2}) in \cref{p.properties}), 
we obtain
\begin{equation}
\cP = \cQ^{-1/2} (\cQ^{1/2} \cR \cQ^{1/2})^{1/2} \cQ^{-1/2} \, .
\end{equation}
One checks directly that this formula solves the invariance equation \eqref{e.invariance},
and thus gives the unique positive solution.

To estimate $\cP$, we follow the steps of \cite{PT}.
By the first inequality in \eqref{e.almost_inv_smart}
and property~(\ref{i.prop1}) in \cref{p.properties},
we have $\cQ^{1/2} \cR \cQ^{1/2} < e^{2 \epsilon} \cQ^2$.
So, by property~(\ref{i.prop3}) in that \lcnamecref{p.properties},  
$(\cQ^{1/2} \cR \cQ^{1/2})^{1/2} < e^\epsilon \cQ$.
Applying property~(\ref{i.prop1}) again, 
we obtain $\cP < e^\epsilon \cI$,
where $\cI$ is the identity automorphism.
This means that $\|P(\omega)\|<e^\epsilon$ for every $\omega$.
An analogous argument starting from the second inequality in \eqref{e.almost_inv_smart}
gives $\m(P(\omega)) > e^{-\epsilon}$ for every $\omega$.
This shows that $\cP$ is close to the identity,
and therefore the automorphism $\tilde \cA := \cA \cP$ is close to $\cA$.
As we have seen, $\tilde \cA$ preserves the new Riemannian metric,
thus completing the proof of the theorem.
\end{proof}

\begin{remark}\label{r.non_positive}
Equation~\eqref{e.invariance} obviously has infinitely many solutions $\cP$, 
not all of them close to the identity.
As we have seen, restricting to positive automorphisms we have a unique solution, 
which is close to the identity and varies continuously with the data. 

In \cite{BN}, we obtain a generalization of \cref{t.zero}
to cocycles of isometries of symmetric spaces of non-positive curvature.
If specialized to the present situation, the construction presented in \cite{BN}
is the same as the one given here for part~(\ref{i.z2}), thus ``explaining''
the efficiency of positive matrices.
\closeremark 
\end{remark}

\begin{remark}\label{r.family}
Notice that the Riemannian norm and the perturbed automorphism
constructed in the proof of \cref{t.zero} 
depend continuously on the parameter $\epsilon$
and also on the automorphism $\cA$ itself.
These properties are relevant for the applications obtained in \cite{ABD2}.
\closeremark 
\end{remark}

\subsection{Conjugacy} \label{ss.conjugacy}

Let us put \cref{t.zero} under a different perspective.

\medskip  

Two automorphisms $\cA$, $\cB \in \Aut(V,F)$ are said to be 
\emph{conjugate} 
if there exists $\cU \in \Aut(V,\id)$ such that $\cA = \cU \cB \cU^{-1}$.
(In the case of a trivial vector bundle, we say that the two linear cocycles are 
\emph{cohomologous}.)

Fixed a Riemannian metric on $V$, we say that an automorphism $\cA$ is \emph{isometric} 
if it preserves this metric.
(In the case of a trivial vector bundle, the cocycle will take values in the orthogonal group,
i.e., it will be a \emph{rotation cocycle}.)

\medskip

Then we have:

\begin{theorem}\label{t.corol}
Fix a Riemannian metric on the vector bundle $V$.
Assume that $\cA \in \Aut(V,F)$ has uniform subexponential growth.
Then:
\begin{enumerate}
\item\label{i.c1}
There exists an automorphism 
conjugate to $\cA$ that is 
close to an isometric automorphism. 
More precisely, every neighborhood of the set of isometric automorphisms contains a conjugate 
of $\cA$. 


\item\label{i.c2} 
There exists an automorphism 
close to $\cA$ 
that is conjugate to an isometric automorphism. More precisely, every neighborhood of $\cA$ contains 
a conjugate of an isometric automorphism.
\end{enumerate}
\end{theorem}


For $\SL(2,\R)$-cocycles
and under extra assumptions on the dynamics $F$,
the result above was shown by Avila, Bochi and Damanik 
as a step in the proofs of their results about spectra of Schr\"odinger operators, see
\cite{ABD1, ABD2}.\footnote{However, they haven't explicitly stated the result: see the 
proof of Theorem~1 in \cite{ABD1} and Proposition~6.3 in \cite{ABD2}.}

In the case of cocycles (i.e., trivial vector bundles), it is natural to
look for conditions under which we can improve the conclusion of \cref{t.corol}(\ref{i.c2})
and find a perturbed cocycle cohomologous to a constant rotation, or even to the identity.
The case of $\SL(2,\R)$-cocycles is studied in \cite{ABD2}.

\begin{proof}[Proof of \cref{t.corol}]
Let $\epsilon >0$ be small.
We follow the notation of the proof of \cref{t.zero}.

It follows from \eqref{e.def_R} that 
$\tribar{v}_\omega = \|R(\omega)^{1/2} v\|_\omega$ for every $v\in V_\omega$.
Let $\cB := \cR^{1/2} \cA \cR^{-1/2}$.
Then, by \eqref{e.almost_inv_metric},
$$
e^{-\epsilon} \|v\|_\omega < \|B(\omega) v \|_{F \omega} < e^\epsilon \|v\|_\omega \, .
$$
This implies that $\cB$ is close to an isometric isomorphism,
thus proving part~(\ref{i.c1}).

To prove part~(\ref{i.c2}), it suffices to notice that $\cR^{-1/2} \tilde \cA \cR^{1/2}$
is an isometric isomorphism.
\end{proof}

There is another property which is closely related to what we have seen so far.
Let us say that $\cA \in \Aut(V,F)$ 
is \emph{product-bounded} 
if 
$$
0 <
\inf_{\substack{\omega \in \Omega \\ n \in \Z}} \m(A^n (\omega))
\le 
\sup_{\substack{\omega \in \Omega \\ n \in \Z}} \| A^n (\omega) \| 
< \infty \, .
$$
If an automorphism $\cA$ is conjugate to an isometric automorphism
then $\cA$ is product-bounded, as it is easy to check.
Although product-bounded cocycles are not always conjugate to isometric automorphisms\footnote{See e.g.\ \cite[Exercise 2.9.2]{KatokH}, \cite{Quas}.},
this happens whenever $F$ is minimal,
according to a result shown by Coronel, Navas and Ponce in \cite{CNP1}.\footnote{In the non-minimal case, one can still ensure the existence of a bounded and measurable conjugacy.}

\section{Conformality properties}

\subsection{Extensions of the previous results for the case of coinciding Lyapunov exponents}

The following is an immediate consequence of \cref{t.zero}:

\begin{corollary}\label{c.equal}
Let $\cA \in \Aut(V,F)$ be such that $\lambda^+(\cA) = \lambda^-(\cA) =: \lambda$.
Then there exist an arbitrarily small perturbation $\tilde \cA$ of $\cA$ 
and a Riemannian norm $\tribar{\mathord{\cdot}}$ on $V$ such that
$$
\tribar{\tilde A(\omega) v }_{F \omega} = e^\lambda \tribar{v}_\omega \, , \quad
\text{for all $\omega \in \Omega$, $v\in V_\omega$.}
$$
\end{corollary}

In other words, if all Lyapunov exponents of an automorphism are equal to some $\lambda$, 
then we can perturb it to become conformal with respect to a new Riemannian 
metric; moreover it dilates the metric by the constant factor $e^\lambda$.

Actually, a weaker assumption is sufficient to obtain conformality:

\begin{corollary}\label{c.not_so_equal}
Let $\cA \in \Aut(V,F)$ be such that $\lambda^+(\cA,\mu) = \lambda^-(\cA,\mu)$
for every ergodic probability measure $\mu$ for $F$.
Then there exist an arbitrarily small perturbation $\tilde \cA$ of $\cA$, a Riemannian norm $\tribar{\mathord{\cdot}}$ on $V$,
and a continuous function $\lambda \colon \Omega \to \R$ such that
$$
\tribar{\tilde A(\omega) v }_{F \omega} = e^{\lambda(\omega)} \tribar{v}_\omega \, , \quad
\text{for all $\omega \in \Omega$, $v\in V_\omega$.}
$$
\end{corollary}

See \cite{KalS} for a non-perturbative result with a similar conclusion.

\begin{proof}[Proof of \cref{c.not_so_equal}]
First of all, notice that for any $\cA \in \Aut(V,F)$, we have
$$
\m \left( A(\omega) \right)^d \le \left| \det A(\omega) \right| \le \|A(\omega)\|^d \, ,
$$
where $d$ be the fiber dimension of $V$.
So, by submultiplicativity of norms, 
\begin{equation}\label{e.lambda_sandwich}
\lambda^-(\cA,\mu) \le \int \lambda \, d\mu \le \lambda^+(\cA,\mu) \quad
\text{for every ergodic measure $\mu$,}
\end{equation}
where
\begin{equation}\label{e.lambda_formula}
\lambda(\omega):= \frac{1}{d} \log \left|\det A_i(\omega)\right|.
\end{equation}

Now assume that equalities hold in \eqref{e.lambda_sandwich}.
Let $\cB = e^{-\lambda} \cA$.
Then $\lambda^{\pm}(\cB,\mu) = \lambda^{\pm}(\cA,\mu) -  \int \lambda \, d\mu = 0$ 
for every ergodic measure $\mu$.
By the ``variational principle'' \eqref{e.variational}, this implies 
$\lambda^+(\cB) = \lambda^-(\cB) = 0$.
Therefore, by \cref{t.zero}, there is
a Riemannian norm $\tribar{\mathord{\cdot}}$ on $V$ 
that is preserved by a perturbation $\tilde \cB$ of $\cB$.

Let $\tilde \cA = e^{\lambda} \tilde \cB$.
This is a perturbation of $\cA$ with the desired conformality property.
\end{proof}


\subsection{Existence of conformal subbundles}\label{ss.subbundles}

Using \cref{c.equal} and a theorem from \cite{BV},
we will obtain the following result:

\begin{theorem}\label{t.conformal_ue}
Assume that $F \colon \Omega \to \Omega$ is a uniquely ergodic homeomorphism
with an invariant probability measure of full support.
Then for every automorphism $\cA$ in a dense subset of $\Aut(V,F)$, 
there exist:
\begin{itemize}
\item a Riemannian norm $\tribar{\mathord{\cdot}}$ on $V$;
\item a continuous $\cA$-invariant splitting $V = V^1 \oplus \cdots \oplus V^k$
which is orthogonal with respect to the Riemannian norm;
\item and constants $\lambda_1 > \cdots > \lambda_k$;
\end{itemize}
such that
$$
\tribar{A(\omega)  v_i}_{F \omega} = e^{\lambda_i} \tribar{v_i}_\omega  \, , 
\quad \text{for all $\omega\in \Omega$, $i=1,\dots,k$, $v_i \in V^i_\omega$}.
$$
\end{theorem}

Weakening the assumption of unique ergodicity to minimality,
we have the following result:

\begin{theorem}\label{t.conformal_minimal}
Assume that $F \colon \Omega \to \Omega$ is a minimal homeomorphism
of a compact space of finite dimension\footnote{We say that $\Omega$ has \emph{finite dimension} if it is homeomorphic to a subset of an euclidean space $\R^d$.}.
Then for every $\cA$ in a dense subset of $\Aut(V,F)$,
there exist:
\begin{itemize}
\item a Riemannian norm $\tribar{\mathord{\cdot}}$ on $V$;
\item a continuous $\cA$-invariant splitting $V = V^1 \oplus \cdots \oplus V^k$
which is orthogonal with respect to the Riemannian norm;
\item and continuous  functions $\lambda_1 > \cdots > \lambda_k$ on $\Omega$; 
\end{itemize}
such that
$$
\tribar{A(\omega)  v_i}_{F \omega} = e^{\lambda_i(\omega)} \tribar{v_i}_\omega  \, ,
\quad \text{for all $\omega\in \Omega$, $i=1,\dots,k$, $v_i \in V^i_\omega$}.
$$
\end{theorem}

As we will see, this has a similar proof as \cref{t.conformal_ue},
basically replacing \cref{c.equal} by \cref{c.not_so_equal}  
and the result from \cite{BV} by the result from \cite{Bo}.

\medskip

We expect that \cref{t.conformal_ue,t.conformal_minimal} 
will be useful to answer the following question:
\emph{When can a linear cocycle over a uniquely ergodic or minimal base dynamics
be approximated by a cocycle with a dominated (non-trivial) splitting?}
Results on the $2$-dimensional case were obtained in \cite{ABD1,ABD2}.


\subsection{Proofs}

\begin{proof}[Proof of \cref{t.conformal_ue}] 
Assume that $F \colon \Omega \to \Omega$ 
has a unique invariant probability $\mu$, and its support is $\Omega$.
Take any $\cA \in \Aut(V,F)$;
we will explain how to perturb it so that it has the desired properties.
First, by \cite{BV}, one can perturb $\cA$ so that along $\mu$-almost every orbit,
the Oseledets splitting is trivial or dominated.
Let 
$$
V^1_\omega \oplus \dots \oplus V^k_\omega = V_\omega \, ,  \quad \omega \in \Omega, 
$$
be the \emph{finest dominated splitting} of the cocycle,
that is, the unique
everywhere defined global dominated splitting with a maximal number $k$
of bundles (with $k=1$ if there is no dominated splitting).\footnote{See \cite{BDV} for details on finest dominated splittings.}

We claim that for almost every point, 
there are exactly $k$ different Lyapunov exponents. Indeed, on the one hand,  
there are at least $k$ different exponents because there is a dominated splitting with $k$ bundles.
On the other hand, if there is a positive measure set of points with more than $k$ different
Lyapunov exponents, then select an orbit along which the Oseledets splitting is dominated.
This orbit is dense on $\Omega$ (because the invariant measure has full support).
Since dominated splittings extend to the closure (see \cite{BDV}), 
one gets a global dominated splitting
with more than $k$ bundles; this is a contradiction.

\medskip

For each $i=1,\ldots,k$, let $\cA_i$ be the restriction of $\cA$ to the bundle $V^i$; 
this is a (continuous) vector bundle automorphism. 
By the claim above,
$$
\lambda^+(\cA_i) = \lambda^+(\cA_i, \mu) = \lambda^-(\cA_i, \mu) = \lambda^-(\cA_i) =: \lambda_i \, .
$$
Therefore, by \cref{c.equal}, for each $i$ 
there is a perturbation $\tilde \cA_i$ of $\cA_i$ 
and a Riemannian norm $\tribar{\mathord{\cdot}}_i$ on $V^i$ such that
$$
\tribar{\tilde A_i(\omega) v_i}_{i, F \omega} = e^{\lambda} \tribar{v_i}_{i, \omega} \, ,
\quad \text{for all $\omega \in \Omega$, $v_i \in V^i_\omega$.}
$$

Let $\tribar{\mathord{\cdot}}$ be the Riemannian norm that makes the subbundles orthogonal 
and that coincides with $\tribar{\mathord{\cdot}}_i$ on $V^i$.
Let $\tilde \cA$ be the automorphism of $V$ whose restriction to the subbundles $V^i$
are the automorphisms $\tilde \cA_i$.
This automorphism has the desired properties, thus completing the proof.
\end{proof}

For the proof of \cref{t.conformal_minimal}, we need the following result:

\begin{theorem}\label{t.multidim_AB}
Assume that $F \colon \Omega \to \Omega$ is a minimal homeomorphism of a compact space of finite dimension.
Then every $\cA$ in a residual subset of $\Aut(V,F)$ has the following property: 
the Oseledets splitting with respect to any invariant probability measure
coincides almost everywhere with the finest dominated splitting of $\cA$.
\end{theorem}

This result is proved in full generality in \cite{Bo}. (The case of $\SL(2,\R)$-cocycles was previously considered in
\cite{AB}.) Notice that, as we have seen in the proof of \cref{t.conformal_ue} above, under 
the additional assumption of unique ergodicity, \cref{t.multidim_AB} follows from \cite{BV}. 

\begin{proof}[Proof of \cref{t.conformal_minimal}] 
Assume that $F \colon \Omega \to \Omega$ is minimal.
Take any $\cA \in \Aut(V,F)$;
we will explain how to perturb it so that it has the desired properties.
First, perturb $\cA$ so that it has the property from \cref{t.multidim_AB}.
This means that if 
$$
V^1_\omega \oplus \dots \oplus V^k_\omega = V_\omega \, ,  \quad (\omega \in \Omega)
$$
is the finest dominated splitting of the cocycle
and $\cA_i$ is the restriction of $\cA$ to the bundle $V^i$ then
\begin{equation}\label{e.conclusion_AB}
\lambda^+(\cA_i, \mu) = \lambda^-(\cA_i, \mu) 
\quad \text{for every ergodic probability $\mu$ for~$F$.}
\end{equation}

By \cite{Go}, we can choose a Riemannian metric on $V$ that is \emph{adapted}
to the dominated splitting, which means that
$$
\inf_{\omega \in \Omega}
\frac{\m(A_i(\omega))}{\| A_{i+1}(\omega) \|} > 1,
\quad \text{for every $i=1,2,\dots,k-1$.}
$$
Let $d_i$ be the fiber dimension of $V^i$, and let 
\begin{equation}\label{e.lambda_i_formula}
\lambda_i(\omega):= \frac{1}{d_i} \log \left|\det A_i(\omega)\right|;
\end{equation}
here determinants are computed with respect to the adapted metric,
and in particular
$\lambda_1 > \lambda_2 > \cdots > \lambda_d$ pointwise.\footnote{One can 
avoid using adapted metrics in the proof of \cref{t.conformal_minimal} by using the following fact: if the functions $\lambda_1$, \dots, $\lambda_k$ satisfy $\int \lambda_1 \, d\mu > \dots > \int \lambda_k \, d\mu $ for every ergodic probability $\mu$, then there are functions $\hat{\lambda}_i$ cohomologous to the $\lambda_i$'s such that $\hat\lambda_1 > \cdots > \hat\lambda_d$ pointwise.}

For each $i$, property \eqref{e.conclusion_AB} permits 
us to apply \cref{c.not_so_equal} and find
a perturbation $\tilde \cA_i$ of $\cA_i$
that is conformal with respect to some Riemannian norm $\tribar{\mathord{\cdot}}_i$ on $V^i$.
Recalling formula \eqref{e.lambda_formula} from the proof of \cref{c.not_so_equal},
we see that 
$\tribar{A(\omega)  v}_{i, F \omega} = e^{\lambda_i(\omega)} \tribar{v}_{i,\omega}$
where the function $\lambda_i$ is given by \eqref{e.lambda_i_formula}.


Let $\tribar{\mathord{\cdot}}$ be the Riemannian norm that makes the subbundles orthogonal 
and that coincides with $\tribar{\mathord{\cdot}}_i$ on $V^i$.
Let $\tilde \cA$ be the automorphism of $V$ whose restrictions to the subbundles $V^i$
are the automorphisms $\tilde \cA_i$. 
This automorphism has the desired properties, thus completing the proof.
\end{proof}

\begin{ack}
We thank the referee for several suggestions and corrections.
\end{ack}


\begin{small}

\phantomsection 

\vspace{0.5cm}

\noindent
\begin{minipage}[t]{.4\linewidth}
Jairo Bochi

PUC--Rio

Rua Marqu\^es de S.~Vicente, 225

Rio de Janeiro, Brazil

jairo@mat.puc-rio.br

\href{http://www.mat.puc-rio.br/~jairo}{www.mat.puc-rio.br/$\sim$jairo}
\end{minipage}
\hspace{.2\linewidth}
\begin{minipage}[t]{.4\linewidth}
Andr\'es Navas

Universidad de Santiago

Alameda 3363, Estaci\'on Central

Santiago, Chile

andres.navas@usach.cl
\end{minipage}

\end{small}

\end{document}